\documentclass[a4paper,11pt,reqno]{smfart}
\usepackage{amssymb,amsmath,mathrsfs,graphics,graphicx,mathptm,float,lscape}
\usepackage[T1]{fontenc}
\usepackage[english, francais]{babel}
\usepackage[all]{xy}
\theoremstyle{plain}
\newtheorem{thm}{Theorem}[section]

\newtheorem{cor}[thm]{Corollary}
\newtheorem{theoalph}{Theorem}

\theoremstyle{definition}
\newtheorem*{defi}{Definition}

\newtheorem{rem}[thm]{Remark}

\usepackage[colorlinks, linktocpage, 
citecolor = blue, linkcolor = blue]{hyperref}

\def\og{\leavevmode\raise.3ex\hbox{$\scriptscriptstyle\langle\!\langle$~}}
\def\fg{\leavevmode\raise.3ex\hbox{~$\!\scriptscriptstyle\,\rangle\!\rangle$}}

\addtocounter{section}{0}             
\numberwithin{equation}{section}       

\begin{document}
\selectlanguage{english}
\title{\textsc{Jonqui\`eres} maps and $\mathrm{SL}(2;\mathbb{C})$-cocycles}

\author{Julie \textsc{D\'eserti}}
\address{Institut de Math\'ematiques de Jussieu-Paris Rive Gauche, UMR $7586$, Universit\'e Paris Diderot, B\^atiment Sophie Germain, Case $7012$, $75205$ Paris Cedex $13$, France.}
\email{deserti@math.univ-paris-diderot.fr}

\maketitle

\begin{abstract}
We start the study of the family of birational maps $(f_{\alpha,\beta})$ of $\mathbb{P}^2_\mathbb{C}$ in~\cite{Deserti}. For "$(\alpha,\beta)$ well chosen" of modulus $1$ the centraliser of $f_{\alpha,\beta}$ is trivial, the topological entropy of~$f_{\alpha,\beta}$ is $0$, there exist two domains of linearisation: in the first one the closure of the orbit of a point is a torus, in the other one the closure of the orbit of a point is the union of two circles. On~$\mathbb{P}^1_\mathbb{C}\times \mathbb{P}^1_\mathbb{C}$ any $f_{\alpha,\beta}$ can be viewed as a cocyle; using recent results about $\mathrm{SL}(2;\mathbb{C})$-cocycles~(\cite{Avila}) we determine the \textsc{Lyapunov} exponent of the cocyle associated to $f_{\alpha,\beta}$. 
\\

\noindent\emph{$2010$ Mathematics Subject Classification. --- $37$F$10$, $14$E$07$}
\end{abstract}

\section*{Introduction}

In this article we deal with a family of birational maps $(f_{\alpha,\beta})$ given by 
\[
f_{\alpha,\beta}\colon\mathbb{P}^2_\mathbb{C}\dashrightarrow\mathbb{P}^2_\mathbb{C}\quad\quad\quad (x:y:z)\dashrightarrow\big((\alpha x+y)z:\beta y(x+z):z(x+z)\big)
\]
where $\alpha$, $\beta$ denote two complex numbers with modulus $1$, case where we know almost nothing about the dynamics. Let us consider the set $\Omega$ of pairs of complex numbers of modulus $1$ that satisfy diophantine condition. The family $(f_{\alpha,\beta})$ satisfies the following properties (\cite{Deserti}):
\begin{itemize}
\item[$\bullet$] for $(\alpha,\beta)\in\Omega$ the centraliser of $f_{\alpha,\beta}$, that is the set of birational maps of $\mathbb{P}^2_\mathbb{C}$ that commutes with $f_{\alpha,\beta}$, is isomorphic to $\mathbb{Z}$;

\item[$\bullet$] the topological entropy of $f_{\alpha,\beta}$ is $0$;

\item[$\bullet$] rotation domains of ranks $1$ and $2$ coexist: there is a domain of linearisation where the orbit of a generic point under $f_{\alpha,\beta}$ is a torus, and there is an other domain of linearisation where the orbit of a generic point under $f_{\alpha,\beta}^2$ is a circle.
\end{itemize}

\medskip

We can also see $f_{\alpha,\beta}$ on $\mathbb{P}^1_\mathbb{C}\times\mathbb{P}^1_\mathbb{C}$  $\big($since all the computations of \cite{Deserti} have been done in an affine chart they may all be carried on $\mathbb{P}^1_\mathbb{C}\times\mathbb{P}^1_\mathbb{C}$$\big)$; the sets~$\mathbb{P}^1_\mathbb{C}\times\mathbb{S}^1_\rho$, where $\mathbb{S}^1_\rho=\{y\in\mathbb{C}\,\vert\,\vert y\vert=\rho\}$, are invariant. 

\smallskip

Let us define $A_n^{\alpha,\rho}\colon\mathbb{S}^1_\rho\rightarrow\mathrm{M}(2;\mathbb{C})$ given in terms of $A^{\alpha,\rho}(y)=\left[\begin{array}{cc} \alpha & y \\ 1 & 1 \end{array}\right]$ by 
\[
A_n^{\alpha,\rho}(\cdot)=A^{\alpha,\rho}(\beta^n\,\cdot)A^{\alpha,\rho}(\beta^{n-1}\,\cdot)\ldots A^{\alpha,\rho}(\beta\, \cdot)A^{\alpha,\rho}(\cdot).
\]

To compute $f_{\alpha,\beta}^n(x,y)$ is equivalent to compute~$A_n^{\alpha,\rho}(y)$ as soon as $f^k_{\alpha,\beta}(x,y)\not=(-1,\alpha)$ for any $1\leq k\leq n$.

\medskip

Using \cite{Avila} we are able to determine the \textsc{Lyapunov} exponent of the cocycle~$(A^{\alpha,\rho},\beta)$:

\begin{theoalph}\label{Thm:main}
{\sl The \textsc{Lyapunov} exponent of $(A^{\alpha,\rho},\beta)$ is 
\begin{itemize}
\item[$\bullet$] positive as soon as $\rho>1$;

\item[$\bullet$] zero as soon as $\rho\leq 1$.
\end{itemize}

\smallskip

More precisely $f_{\alpha,\beta}$ is semi conjugate to $\left(\frac{\alpha x+y^2}{x+1},\beta^{1/2}y\right)$ and the \textsc{Lyapunov} exponent of the cocycle $\big(B^{\alpha,\rho},\beta^{1/2}\big)$, where 
\[
B^{\alpha,\rho}(y)=\left[\begin{array}{cc} \alpha & y^2 \\ 1 & 1 \end{array}\right],
\]
is equal to $\max(0,\ln\rho)$.}
\end{theoalph}

\bigskip

In the next section we introduce the family $(f_{\alpha,\beta})$ and its properties (\S\ref{Sec:falphabeta}). Then we deal with the recent works of \textsc{Avila} on $\mathrm{SL}(2;\mathbb{C})$-cocyles. In the last section we give the proof of Theorem \ref{Thm:main} (\emph{see} \S\ref{Sec:cocyles}). Let us explain the sketch of it. We associate to $\big(B^{\alpha,\rho},\beta^{1/2}\big)$ a cocycle $\big(\widetilde{B}^{\alpha,\rho},\beta^{1/2}\big)$ that belongs to~$\mathrm{SL}(2;\mathbb{C})$. We first determine 
\[
\displaystyle\lim_{\rho\to 0}L\big(\widetilde{B}^{\alpha,\rho},\beta^{1/2}\big),
\]
and then 
\[
\displaystyle\lim_{\rho\to +\infty}L\big(\widetilde{B}^{\alpha,\rho},\beta^{1/2}\big)
\]
where $L(C,\gamma)$ denotes the \textsc{Lyapunov} exponent of the $\mathrm{SL}(2;\mathbb{C})$-cocyle $(C,\gamma)$. In both cases, we get $0$. Using \cite[Theorem 5]{Avila} we obtain that $L\big(\widetilde{B}^{\alpha,\rho},\beta^{1/2}\big)$ vanishes everywhere; it allows us to determine $L\big(A^{\alpha,\rho},\beta\big)$ since 
\[
L\big(B^{\alpha,\rho}(y),\beta^{1/2}\big)=L\big(\widetilde{B}^{\alpha,\rho}(y),\beta^{1/2}\big)+\max(0,\ln\rho),
\]
and since $\big(A^{\alpha,\rho},\beta\big)$ and $\big(\beta^{1/2},B^{\alpha,\rho}\big)$ are conjugate.

\medskip

\subsection*{Acknowledgment} I would like to thank Artur \textsc{Avila} for very helpful discussions, Dominique \textsc{Cerveau} for his constant support, and Serge \textsc{Cantat} for his remarks. Thanks also to the referee whose comments help me to improve the text.

\section{Some properties of the family $(f_{\alpha,\beta})$}\label{Sec:falphabeta}

A \textbf{\textit{rational map}} $\phi$ from $\mathbb{P}^2_\mathbb{C}$ into itself is a map of the form 
\[
(x:y:z)\dashrightarrow\big(\phi_0(x,y,z):\phi_1(x,y,z):\phi_2(x,y,z)\big),
\]
where the $\phi_i$'s are some homogeneous polynomials of the same degree without common factor;~$\phi$ is \textbf{\textit{birational}} if it admits an inverse of the same type. We will denote by $\mathrm{Bir}(\mathbb{P}^2_\mathbb{C})$ the group of birational maps of $\mathbb{P}^2_\mathbb{C}$, also called the \textbf{\textit{\textsc{Cremona} group}}. The \textbf{\textit{degree}} of $\phi$, denoted $\deg\phi$, is the degree of the $\phi_i$'s. The degree is not a birational invariant: $\deg\psi\phi\psi^{-1}\not=\deg\phi$ for generic birational maps $\phi$ and $\psi$. The \textbf{\textit{first dynamical degree}} of $\phi$ given by 
\[
\lambda(\phi)=\lim_{n\rightarrow +\infty}\big(\deg \phi^n\big)^{1/n},
\]
is a birational invariant; it is strongly related to the topological 
entropy $h_{\text{top}}(\phi)$ of $\phi$ (\emph{see} \cite{Gromov, Yomdin}) 
\begin{equation}\label{GromovYomdin}
h_{\text{top}}(\phi)\leq\log\lambda(\phi)
\end{equation}

Any birational map $\phi$ admits a resolution 
\[
\xymatrix{& \mathrm{S}\ar[rd]^{\pi_2}\ar[ld]_{\pi_1}&\\
\mathbb{P}^2_\mathbb{C}\ar@{-->}[rr]_{\phi}&&\mathbb{P}^2_\mathbb{C},
}
\]
where $\pi_1$, $\pi_2\colon\mathrm{S}\to\mathbb{P}^2_\mathbb{C}$ are sequences of blow-ups (\emph{see} \cite{Beauville} for example). The resolution is \textbf{\textit{minimal}} if and only if no $(-1)$-curve of $\mathrm{S}$ is contracted by both $\pi_1$ and $\pi_2$. The \textbf{\textit{base-points}} of $\phi$ are the points blown-up in~$\pi_1$, which can be points of $\mathbb{P}^2_\mathbb{C}$ or infinitely near points. We denote by $\mathfrak{b}(\phi)$ the number of such points, which is also equal to the difference of the ranks of $\mathrm{Pic}(\mathrm{S})$ and $\mathrm{Pic}(\mathbb{P}^2_\mathbb{C})$, and thus equals to~$\mathfrak{b}(\phi^{-1})$. The \textbf{\textit{dynamical number of base-points of $\phi$}} introduced in~\cite{BlancDeserti} is by definition 
\[
\mu(\phi)=\displaystyle\lim_{n\rightarrow +\infty} \frac{\mathfrak{b}(\phi^n)}{n};
\]
it is a real positive number that satisfies $\mu(\phi^n)=\vert n\, \mu(\phi)\vert$ for any $n\in \mathbb{Z}$, $\mu(\psi\phi\psi^{-1})=\mu(\phi)$, and allows us to give a characterization of birational maps conjugate to automorphisms:

\begin{thm}[\cite{BlancDeserti}]\label{thm:caraut}
{\sl Let $\mathrm{S}$ be a smooth projective surface; the birational map $\phi\in\mathrm{Bir}(\mathrm{S})$ is conjugate to an automorphism of a smooth projective surface if and only if $\mu(\phi)=0$.}
\end{thm}

The behavior of $\phi\in\mathrm{Bir}(\mathbb{P}^2_\mathbb{C})$ is strongly related to the behavior of~$\big(\deg\phi^n\big)_{n\in\mathbb{N}}$ (\emph{see} \cite{Gizatullin, DillerFavre, BlancDeserti}); up to birational conjugacy exactly one of the following holds:
\begin{enumerate}
\item the sequence $\big(\deg\phi^n\big)_{n\in\mathbb{N}}$ is bounded and either $\phi$ is of finite order, or $\phi$ is an automorphism of~$\mathbb{P}^2_\mathbb{C}$;

\item there exists an integer $k$ such that 
\[
\displaystyle\lim_{n\rightarrow +\infty} \frac{\deg \phi^n}{n}=k^2\,\frac{\mu(\phi)}{2}
\]
and $\phi$ is not an automorphism;

\item there exists an integer $k\geq 3$ such that 
\[
\displaystyle\lim_{n\rightarrow +\infty} \frac{\deg \phi^n}{n^2}=k^2\,\frac{\kappa(\phi)}{9}
\]
where $\kappa(\phi)\in \mathbb{Q}$ is a birational invariant, and $\phi$ is an automorphism;

\item the sequence $\big(\deg \phi^n\big)_{n\in\mathbb{N}}$ grows exponentially (\emph{see} \cite{DillerFavre} for more precise dynamical pro\-perties).
\end{enumerate}

In the first three cases $\lambda(\phi)=1$, in the last one $\lambda(\phi)>1$. 
In case 2. (resp. 3.) the map~$\phi$ preserves a unique fibration which is rational (resp. elliptic).

\medskip

In case 1. (resp. 2., resp. 3, resp. 4) we say that $\phi$ is \textbf{\textit{elliptic}} (resp. a \textbf{\textit{\textsc{Jonqui\`eres} twist}}, resp. an \textbf{\textit{\textsc{Halphen} twist}}, resp. \textbf{\textit{hyperbolic}}).

\smallskip

Let us give some examples. Let 
\[
\phi(x,y)=\left(\frac{a(y)x+b(y)}{c(y)x+d(y)},\frac{\alpha y+\beta}{\gamma y+\delta}\right)
\]
be an element of the \textbf{\textit{\textsc{Jonqui\`eres} group}} $\mathrm{PGL}(2;\mathbb{C}(y))\rtimes\mathrm{PGL}(2;\mathbb{C})$; either $\phi$ is elliptic (for instance $\phi\colon(x:y:z)\dashrightarrow(yz:xz:xy)$), or $\phi$ is a \textsc{Jonqui\`eres} twist (for example $\phi\colon(x:y:z)\dashrightarrow(xz:xy:z^2)$ for which the unique invariant fibration is $y/z=$ constant). The map  
\[
\phi\colon\mathbb{P}^2_\mathbb{C}\dashrightarrow\mathbb{P}^2_\mathbb{C}\quad\quad\quad (x:y:z)\dashrightarrow\big((2y+z)(y+z):x(2y-z):2z(y+z)\big)
\] 
is an \textsc{Halphen} twist (\cite[Proposition 9.5]{DillerFavre}). \textsc{H\'enon} automorphisms give by homogeneization examples of hyperbolic maps.

\bigskip

Clearly elliptic birational maps have a poor dynamical behavior contrary to hyperbolic ones. The study of automorphisms of positive entropy is strongly related with birational maps of $\mathbb{P}^2_\mathbb{C}$:

\begin{thm}[\cite{Cantat}]
{\sl Let $\mathrm{S}$ be a compact complex surface that carries an automorphism $\phi$ of positive topological entropy.
\begin{itemize}
\item[$\bullet$] Either the \textsc{Kodaira} dimension of $\mathrm{S}$ is zero and $\phi$ is conjugate to an automorphism on the unique minimal model of $\mathrm{S}$ that necessarily is a torus, or a K$3$ surface or an \textsc{Enriques} surface;

\item[$\bullet$] or the surface $\mathrm{S}$ is a non-minimal rational one, isomorphic to $\mathbb{P}^2_\mathbb{C}$ blown up at $n$ points, $n\geq 10$, and $\phi$ is conjugate to a birational map of $\mathbb{P}^2_\mathbb{C}$.
\end{itemize}
}
\end{thm}

This yields many examples of hyperbolic birational maps for which we can establish a lot of dyna\-mical properties (\cite{McMullen, BedfordKim1, BedfordKim2, BedfordKim3, BedfordKim4, Diller, DesertiGrivaux}). 

\medskip

Another way to measure chaos is to look at the size of centralisers. Let us give two examples. The polynomial automorphisms of $\mathbb{C}^2$ having rich dynamics are \textsc{H\'enon} maps; furthermore a polynomial automorphism of $\mathbb{C}^2$ is a \textsc{H\'enon} one if and only if its centraliser is countable. Let us now consider rational maps on $\mathbb{S}^1$; if the centraliser of such maps is not trivial\footnote{The centraliser of a map $\phi$ is trivial if it coincides with the iterates of $\phi$.}, then the \textsc{Julia} set is "special". The centraliser of an elliptic birational map of infinite order is uncountable (\cite{BlancDeserti}). The centralisers of \textsc{Halphen} twists are described in \cite{Gizatullin}. The centraliser of an hyperbolic map is countable (\cite{Cantat:annals}). In \cite{CerveauDeserti} we end the story by studying centralisers of \textsc{Jonqui\`eres} twists. If the fibration is fiberwise invariant, then the centraliser is uncountable ; but if it isn't, then generically the centraliser is isomorphic to $\mathbb{Z}$. We don't know a lot about dynamics of these maps, in this article we will thus focus on a family of such maps. We consider the \textsc{Jonqui\`eres} maps 
\[
f_{\alpha,\beta}\colon\mathbb{P}^2_\mathbb{C}\dashrightarrow\mathbb{P}^2_\mathbb{C}\quad\quad\quad (x:y:z)\dashrightarrow\big((\alpha x+y)z:\beta y(x+z):z(x+z)\big)
\]
where $\alpha$, $\beta$ denote two complex numbers with modulus $1$. The base-points of $f_{\alpha,\beta}$ are 
\[
(1:0:0),\quad(0:1:0),\quad(-1:\alpha:1).
\]

Any $f_{\alpha,\beta}$ preserves a rational fibration (the fibration $y=$ constant in the affine chart $z=1$). Each element of the fa\-mily~$(f_{\alpha,\beta})$ has first dynamical degree $1$ hence topological entropy zero (\ref{GromovYomdin}); more precisely one has (\cite[Example 4.3]{BlancDeserti}) 
\[
\mu( f_{\alpha,\beta})=\frac{1}{2}
\]
so $f_{\alpha,\beta}$ is not conjugate to an automorphism (Theorem \ref{thm:caraut}). The centralizer of $f_{\alpha,\beta}$ is isomorphic to $\mathbb{Z}$ (\emph{see} \cite[Theorem~1.6]{Deserti}). The idea of the proof is the following: the point $p = (1 : \alpha : 1)$ is blown-up onto a fiber of the fibration $y=$ constant. Let $\psi$ be an element of 
\[
\mathrm{Cent}(f_{\alpha,\beta})=\big\{g\in\mathrm{Bir}(\mathbb{P}^2_\mathbb{C})\,\vert\, g\circ f_{\alpha,\beta}=f_{\alpha,\beta}\circ g\big\};
\]
since~$\psi$ blows down a finite number of curves there exists a positive integer~$k$ (chosen minimal) such that~$f_{\alpha,\beta}^k(p)$ is not blown down by $\psi$. Replacing $\psi$ by $\widetilde{\psi}=\psi f_{\alpha,\beta}^{k-1}$ one gets that $\widetilde{\psi}(p)$ is an indeterminacy point of~$f_{\alpha,\beta}$. In other words $\widetilde{\psi}$ permutes the indeterminacy points of $f_{\alpha,\beta}$. A more precise study allows us to establish that $p$ is fixed by $\widetilde{\psi}$. The pair $(\alpha,\beta)$ being in $\Omega$, the closure of the negative orbit of $p$ under the action of $f_{\alpha,\beta}$ is \textsc{Zariski} dense; since $\widetilde{\psi}$ fixes any element of the orbit of $p$ one obtains $\widetilde{\psi}=\mathrm{id}$.

\medskip

Let us recall that if $\psi$ is an automorphism on a compact complex manifold $\mathrm{M}$, the \textbf{\textit{\textsc{Fatou} set}}~$\mathcal{F}(\psi)$ of $\psi$ is the set of points that have a neighborhood $\mathcal{V}$ such that $\big\{f^n_{\vert\mathcal{V}}\,\vert\, n\in\mathbb{N}\big\}$ is a normal family. Set
\[
\mathcal{G}(\mathcal{U})=\big\{\phi\colon\mathcal{U}\to\overline{\mathcal{U}}\,\vert\,\phi=\lim_{n_j\to +\infty}\psi^{n_j}\big\};
\]
we say that $\mathcal{U}$ is a \textbf{\textit{rotation domain}} if $\mathcal{G}(\mathcal{U})$ is a subgroup of $\mathrm{Aut}(\mathcal{U})$. An equivalent definition is the following: a component $\mathcal{U}$ of $\mathcal{F}(\psi)$ which is invariant by $\psi$ is a rotation domain if $\psi_{\vert\mathcal{U}}$ is conjugate to a linear rotation. If~$\mathcal{U}$ is a rotation domain, $\mathcal{G}(\mathcal{U})$ is a compact \textsc{Lie} group, and the action of $\mathcal{G}(\mathcal{U})$ on $\mathcal{U}$ is analytic real. Since $\mathcal{G}(\mathcal{U})$ is a compact, infinite, abelian \textsc{Lie} group, the connected component of the identity of $\mathcal{G}(\mathcal{U})$ is a torus of dimension $0\leq d\leq \dim_\mathbb{C}\mathrm{M}$. The integer $d$ is the \textbf{\textit{rank of the rotation domain}}. The rank coincides with the dimension of the closure of a generic orbit of a point in $\mathcal{U}$.

\medskip

We can also see $f_{\alpha,\beta}$ on $\mathbb{P}^1_\mathbb{C}\times\mathbb{P}^1_\mathbb{C}$ and that is what we will do in the sequel $\big($since all the computations of \cite{Deserti} have been done in an affine chart they may all be carried on $\mathbb{P}^1_\mathbb{C}\times\mathbb{P}^1_\mathbb{C}$$\big)$; the sets~$\mathbb{P}^1_\mathbb{C}\times\mathbb{S}^1_\rho$ are invariant. In \cite{Deserti} we show that there are two rotation domains for $f^2_{\alpha,\beta}$, one of rank $1$, and the other one of rank $2$\footnote{There already exists an example of automorphism of positive entropy with rotation domains of rank $1$ and $2$ (\emph{see} \cite{BedfordKim2}), but $f_{\alpha,\beta}$ is not conjugate to an automorphism on a rational surface.}; in the first case we give here a more precise statement than in~\cite{Deserti}:

\begin{thm}\label{Thm:old}
{\sl Assume that $(\alpha,\beta)$ belongs to $\Omega$.

There exists a strictly po\-sitive real number $r$ such that $f_{\alpha,\beta}$ is conjugate to $(\alpha x,\beta y)$ on $\mathbb{P}^1_\mathbb{C}~\times~\mathbb{D}(0,r)$ where $\mathbb{D}(0,r)$ denotes the disk centered at the origin with radius $r$.

There exists a strictly po\-sitive real number $\widetilde{r}$ such that $f_{\alpha,\beta}^2$ is conjugate to $\left(\frac{x}{\beta},\frac{z}{\beta^2}\right)$ on $\mathbb{P}^1_\mathbb{C}\times\mathbb{D}(0,\widetilde{r})$.}
\end{thm}

\begin{rem}
The point $(\alpha-1,0)$ is also a fixed point of $f_{\alpha,\beta}$ where the behavior of~$f_{\alpha,\beta}$ is the same as near $(0,0)$.
\end{rem}

\begin{proof} 
The first assertion is proved in \cite{Deserti}.

Let us consider the map $\psi(x,z)=\left(\frac{a(z)x+b(z)}{c(z)x+1},z\right)$. The equation 
\[
\psi^{-1}f_{\alpha,\beta}^2\psi=\left(\frac{x}{\beta},\frac{z}{\beta^2}\right)
\]
yields 
\begin{eqnarray}\label{eq1}
&&\beta\, a\big(\beta^{-2}\,z\big)c(z)+\beta\, a\big(\beta^{-2}\,z\big)a(z)-c\big(\beta^{-2}\,z\big)a(z)+\alpha\, a\big(\beta^{-2}\,z\big)a(z)\nonumber\\
&&\hspace{0.5cm}+z\big(\alpha^2\,a\big(\beta^{-2}\,z\big)c(z)-\alpha\, c\big(\beta^{-2}\,z\big)c(z)-c\big(\beta^{-2}\,z\big)c(z)-c\big(\beta^{-2}\,z\big)a(z)\big)=0,
\end{eqnarray}
\begin{eqnarray}\label{eq2}
&&\beta\, a\big(\beta^{-2}\,z\big)-\beta\, a(z)
+z\big(\alpha^2\,a\big(\beta^{-2}\,z\big)-\alpha\beta\, c(z)-\beta\, c(z)-\beta\, a(z)-\alpha\, c\big(\beta^{-2}\,z\big)-c\big(\beta^{-2}\,z\big)\big)\nonumber\\
&&\hspace{0.5cm}
+\beta(\alpha+\beta)\, a(z)b\big(\beta^{-2}\,z\big) 
+(\alpha+\beta)\, b(z)a\big(\beta^{-2}\,z\big)+\beta^2\,b\big(\beta^{-2}\,z\big)c(z)-b(z)c\big(\beta^{-2}\,z\big)\nonumber\\
&&\hspace{0.5cm}+z\big(\alpha^2\,\beta b\big(\beta^{-2}\,z\big)c(z)-b(z)c\big(\beta^{-2}\,z\big)\big)=0
\end{eqnarray}
and
\begin{eqnarray}\label{eq3}
(\alpha+1)\, z+ b(z)-\beta\, b\big(\beta^{-2}\,z\big)-\alpha^2 \,zb\big(\beta^{-2}\,z\big)+ z b(z)-(\alpha+\beta)\,b\big(\beta^{-2}\,z\big)b(z)=0
\end{eqnarray}

Let us set 
\[
a(z)=\sum_{i\geq 0}a_iz^i,\quad\quad\quad b(z)=\sum_{i\geq 0}b_iz^i,\quad\quad\quad c(z)=\sum_{i\geq 0}c_iz^i.
\]
We easily get $a_0=1-\beta$, $b_0=0$ and $c_0=\alpha+\beta$.

Relation (\ref{eq3}) implies that 
\[
b_1=\frac{\beta(1+\alpha)}{1-\beta}\quad\quad\quad\&\quad\quad\quad \beta\, b_\nu\,\left(1-\beta^{1-2\nu}\right)+F_i(b_i\,\vert\,i<\nu)=0 \quad\forall\,\nu>1,
\]
(\ref{eq2}) yields 
\begin{eqnarray*}
a_\nu\left(\beta^{1-2\nu}-\beta\right)+b_\nu\left((\alpha+\beta)a_0\Big(1+\beta^{1-2\nu}\Big)+c_0\Big(\beta^{2-2\nu}-1\Big)\right)
+G_i(a_i,\,b_i,\,c_i\,\vert\,i<\nu)=0
\end{eqnarray*}
and (\ref{eq1}) to 
\begin{eqnarray*}
c_\nu a_0\left(\beta-\beta^{-2\nu}\right)+a_\nu\left((\alpha+\beta)a_0\Big(1+\beta^{-2\nu}\Big)+c_0\Big(\beta^{1-2\nu}-1\Big)\right)+H_i(a_i,\,b_i,\,c_i\,\vert\,i<\nu)=0
\end{eqnarray*}
where the $F_i$'s, $G_i$'s and $H_i$'s denote universal polynomials; this allows to compute $b_\nu$, $a_\nu$ and $c_\nu$. Thus we get a formal conjugacy of $f_{\alpha,\beta}^2$ to its linear part. Since this linear part satisfies a R\"ussmann condition (\emph{see} \cite[Theorem 2.1]{Russmann} condition (2)), according to \cite[Theorem 2.1]{Russmann} any formal linearizing map conjugating $f_{\alpha,\beta}^2$ to its linear part is convergent on a polydisc. 
\end{proof}
 
\section{About $\mathrm{SL}(2;\mathbb{C})$-cocycles}\label{Sec:cocyles}
A (one-frequency, analytic) \textbf{\textit{quasiperiodic $\mathrm{SL}(2;\mathbb{C})$-cocycle}} is a pair $(A,\beta)$, where $\beta\in~\mathbb{R}$ and
\[
A\colon\mathbb{S}^1_1\to\mathrm{SL}(2;\mathbb{C})
\]
is analytic, and defines a linear skew product acting on $\mathbb{C}^2\times\mathbb{S}^1_1$ by 
\[
(x,y)\mapsto(A(y)\cdot x,\beta y). 
\]
The iterates of the cocyle are given by $(A_n,n\beta)$ where $A_n$ is given by 
\[
A_n(y)=A\big(\beta^{n-1}y\big)\ldots A(y)\quad n\geq 1,\quad A_0(y)=\mathrm{id},\quad A_{-n}(y)=A_n(\beta^{-n}y)^{-1}.
\]
The \textbf{\textit{\textsc{Lyapunov} exponent}} $L(A,\beta)$ of a quasiperiodic $\mathrm{SL}(2;\mathbb{C})$-cocycle $(A,\beta)$ is given by 
\[
\displaystyle\lim_{n\to +\infty}\frac{1}{n}\int_{\mathbb{S}^1_1}\mathrm{ln}\,\vert\vert A_n(y)\vert\vert\, \mathrm{d}y.
\]
A quasiperiodic $\mathrm{SL}(2;\mathbb{C})$-cocycle $(A,\beta)$ is \textbf{\textit{uniformly hyperbolic}} if there exist ana\-lytic functions 
\[
u,\,s \colon \mathbb{S}^1_1\to\mathbb{P}^2_\mathbb{C},
\] 
called the \textbf{\textit{unstable and stable directions}}, and $n\geq 1$ such that for any $y\in\mathbb{S}^1_1$, 
\[
A(y)\cdot u(y)=u(\beta y)\qquad A(y)\cdot s(y)=s(\beta y), 
\] 
and for any unit vector $x\in s(y)$ (resp. $x\in u(y)$) we have $\vert\vert A_n(y)\cdot x\vert\vert<1$ (resp. $\vert\vert A_n(y)\cdot x\vert\vert>1$). The unstable and stable directions are uniquely characterized by those properties, and clearly $u(y)\not=s(y)$ for any $y\in\mathbb{S}^1_1$. If $(A,\beta)$ is uniformly hyperbolic, then $L(A,\beta)> 0$.
Let us denote by 
\[
\mathcal{U}\mathcal{H}\subset C^\omega\big(\mathrm{SL}(2;\mathbb{C}),\mathbb{S}^1_1\big) 
\] 
the set of $A$ such that $(A,\beta)$ is uniformly hyperbolic. Uniform hyperbolicity is a stable property: $\mathcal{U}\mathcal{H}$ is open, and $A\mapsto L(A,\beta)$ is analytic over $\mathcal{U}\mathcal{H}$ (regularity properties of the \textsc{Lyapunov} exponent are consequence of the regularity of the unstable and stable directions which depend smoothly on both variables). 

\begin{defi}
Let $(A,\beta)$ be a quasiperiodic $\mathrm{SL}(2;\mathbb{C})$-cocycle. If $L(A,\beta)>0$ but $(A,\beta)\not\in\mathcal{U}\mathcal{H}$, then $(A,\beta)$ is \textbf{\textit{nonuniformly hyperbolic}}.
\end{defi}

If $A\in C^\omega\big(\mathrm{SL}(2;\mathbb{C}),\mathbb{S}^1_1\big)$ admits a holomorphic extension to $\vert\mathrm{Im}\, y\vert<\delta$ then for $\vert\varepsilon\vert<\delta$ we can define $A_\varepsilon\in C^\omega\big(\mathrm{SL}(2;\mathbb{C}),\mathbb{S}^1_1\big)$ by
\[
A_\varepsilon(y)=A(y+\mathbf{i}\varepsilon).
\]
The \textsc{Lyapunov} exponent $L(A_\varepsilon,\beta)$ is a convex function of $\varepsilon$. We can thus introduce the following notion. The \textbf{\textit{acceleration}} of a quasiperiodic $\mathrm{SL}(2;\mathbb{C})$-cocyle $(A,\beta)$ is given by 
\[
\omega(A,\beta)=\lim_{\varepsilon\to 0^+}\frac{1}{2\pi\varepsilon}\big(L(A_\varepsilon,\beta)-L(A,\beta)\big).
\]

\begin{rem}\label{rem:decreasing}
The convexity of the \textsc{Lyapunov} exponent in function of $\varepsilon$ implies that the acceleration is decreasing.
\end{rem}

Since the \textsc{Lyapunov} exponent is a convex and continuous function the acceleration is an upper semi-continuous function in $\mathbb{R}\setminus\mathbb{Q}\times C^\omega\big(\mathrm{SL}(2;\mathbb{C}),\mathbb{S}^1_1\big)$. The acceleration is quantized:

\begin{thm}[\cite{Avila}]\label{Thm:acc}
{\sl If $(A,\beta)$ is a $\mathrm{SL}(2;\mathbb{C})$-cocycle with $\beta\in\mathbb{R}\smallsetminus\mathbb{Q}$, then $\omega(A,\beta)$ is always an integer.}
\end{thm}

A direct consequence is the following:

\begin{cor}
{\sl The function $\varepsilon\mapsto L(A_\varepsilon,\beta)$ is a piecewise affine function of $\varepsilon$.}
\end{cor}

It is thus natural to introduce the notion of regularity. A cocycle 
\[
(A,\beta) \in C^\omega\big(\mathrm{SL}(2;\mathbb{C}),\mathbb{S}^1_1\big)\times\mathbb{R}\setminus\mathbb{Q}
\]
is \textbf{\textit{regular}} if $L(A_\varepsilon,\beta)$ is affine for~$\varepsilon$ in a neighborhood of $0$. In other words $(A,\beta)$ is regular if the equality 
\[
L(A_\varepsilon,\beta)-L(A,\beta)=2\pi\varepsilon\omega(A,\beta) 
\]
holds for all~$\varepsilon$ small, and not only for the positive ones. Regularity is equivalent to the acceleration being locally constant near $(A,\beta)$. It is an open condition in  $C^\omega\big(\mathrm{SL}(2;\mathbb{C}),\mathbb{S}^1_1\big)\times\mathbb{R}\setminus\mathbb{Q}$. The following statement gives a characterization of the dynamics of regular cocycles with positive \textsc{Lyapunov} exponent:

\begin{thm}[\cite{Avila}]
{\sl Let $(A,\beta)$ be a $\mathrm{SL}(2;\mathbb{C})$-cocycle with $\beta\in\mathbb{R}\smallsetminus\mathbb{Q}$. Assume that $L(A,\beta)>0$; then $(A,\beta)$ is regular if and only if~$(A,\beta)$ is~$\mathcal{U}\mathcal{H}$.}
\end{thm}

One striking consequence is the following:

\begin{cor}[\cite{Avila}]
{\sl For any $(A,\beta)$ in $C^\omega\big(\mathrm{SL}(2;\mathbb{C}),\mathbb{S}^1_1\big)\times\mathbb{R}\setminus\mathbb{Q}$ there exists $\varepsilon_0$ such that
\smallskip
\begin{itemize}
\item[$\bullet$] $L(A_\varepsilon,\beta)=0$ $($and $\omega(A,\beta)=0)$ for every $0<\varepsilon<\varepsilon_0$, 
\smallskip
\item[$\bullet$] or $(A_\varepsilon,\beta)\in\mathcal{U}\mathcal{H}$ for every $0<\varepsilon<\varepsilon_0$.
\end{itemize}}
\end{cor}

\begin{rem}
Let us mention that there is a link between $\mathrm{SL}(2;\mathbb{C})$-cocycles and \textsc{Schr\"odinger} operators (\emph{see} \cite{Avila} for more details).
\end{rem}

\section{Proof of Theorem \ref{Thm:main}}\label{Sec:main}

Suppose that $\rho\not=1$, and let us consider the cocycle $(B^{\alpha,\rho},\beta^{1/2})$ where 
\[
B^{\alpha,\rho}(y)=\left[\begin{array}{cc} \alpha & y^2 \\ 1 & 1 \end{array}\right].
\]
Since  
\[
\left(\frac{\alpha x+y}{x+1},\beta y\right)(x,y^2)=(x,y^2)\left(\frac{\alpha x+y^2}{x+1},\beta^{1/2}y\right) 
\]
the cocycles $(A^{\alpha,\rho},\beta)$ and $(B^{\alpha,\rho},\beta^{1/2})$ have the same behavior.
Using two different arguments of mono\-dromy (one for $\rho<1$, and the other one for $\rho>1$) we see that there is a continuous determination for the square root of $\det B^{\alpha,\rho}(y)=\alpha-y^2$. Let us set 
\[
\widetilde{B}^{\alpha,\rho}(y)=\frac{1}{\sqrt{\alpha-y^2}}\,B^{\alpha,\rho}(y)\in\mathrm{SL}(2;\mathbb{C})
\]
that is thus defined on two different domains of analyticity.

According to Theorem \ref{Thm:old} one has $L\big(\widetilde{B}^{\alpha,\rho},\beta^{1/2}\big)=0$ when $\rho$ is close to both $0$ and $\infty$.

\medskip

Assume that $L\big(\widetilde{B}^{\alpha,\rho},\beta^{1/2}\big)$ is non constant. When $\widetilde{B}^{\alpha,\rho}$ is holomorphic, so in particular when $\rho<1$ and $\rho>1$, the acceleration is decreasing (Remark \ref{rem:decreasing}); furthermore the acceleration is positive for $\rho<1$ and negative for $\rho>1$ (because $L$ is continuous). Theorem \ref{Thm:acc} thus implies
\[
\omega\big(\widetilde{B}^{\alpha,1^+},\beta^{1/2}\big)-\omega\big(\widetilde{B}^{\alpha,1^-},\beta^{1/2}\big)\leq -2.
\]
By definition of $\widetilde{B}^{\alpha,\rho}$ we have
\begin{eqnarray*}\label{Lyapunov}
L\big(\widetilde{B}^{\alpha,\rho}(y),\beta^{1/2}\big) &=& L\big(B^{\alpha,\rho}(y),\beta^{1/2}\big)-\int_{\mathbb{S}^1_\rho}\ln\sqrt{\alpha-y^2}\,\mathrm{d}y\\
&=&L\big(B^{\alpha,\rho}(y),\beta^{1/2}\big)-\max(0,\ln\rho).
\end{eqnarray*}
Even though $\big(B^{\alpha,\rho}(y),\beta^{1/2}\big)$ is not a $\mathrm{SL}(2;\mathbb{C})$-cocycle, the \textsc{Lyapunov} exponent is still a convex function of $\log\rho$ (\emph{see for example} \cite{AvilaJitomirskayaSadel}). The jump of $\omega(B^{\alpha,\rho}(y),\beta^{1/2})$ is thus $\geq 0$, and the jump for the second term of the right member is $-1$. Therefore the jump of $L\big(\widetilde{B}^{\alpha,\rho}(y),\beta^{1/2}\big)$ is $\geq -1$: contradiction.

\vspace{8mm}

\bibliographystyle{plain}
\bibliography{biblio}
\nocite{}

\end{document}